\documentclass{article}
\usepackage{amssymb,amsmath,amsthm,cite,epic,rotating}
\usepackage[russian]{babel}
\usepackage{pdfsync}
\usepackage[pdftex]{hyperref}
\newtheorem{theorem}{Теорема}
\newtheorem{lemma}{Лемма}
\newtheorem{proposition}{Предложение}
\newtheorem{corollary}{Следствие}
\usepackage[cp1251]{inputenc}
\title{Cyclically regular semigroups}

\author{S. Kublanovsky}

\date{February 11, 2008}

\begin{document}

\maketitle
\begin{center}\textit{TPO ”Severnyi ochag” \\
30, B.Konjushennaja 15 \\
191186, St.Petersburg, Russia \\
E-mail: stas@norths.spb.su}\end{center}
\begin{abstract}
  We study varieties of semigroups related to completely 0-simple semigroup. We present here an algorithmic descriptions of these varieties interms of ”forbidden” semigroups.
  
\end{abstract}

\tableofcontents
\pagebreak

\section{Введение}

Понятие регулярности, введенное  Фон-Нейманом в 1936 г., играет важную роль в теории полугрупп и колец. Здесь мы вводим понятие циклической регулярности. Полугруппу $S$ назовем \emph{циклически регулярной}, если каждый ее циклический элемент регулярен. Элемент $b$ называется \emph{циклическим}, если $b=axa$ для некоторых $a$ из $S$ и $x$ из $S^{1}.$  Потребность такого обобщения вызвана тем обстоятельством, что многообразие, порожденное тем или иным классом регулярных полугрупп, может не состоять из регулярных полугрупп. Однако, при переходе к многообразиям в ряде случаев сохраняется свойство циклической регулярности. Это факт имеет место для  0-простых периодических полугрупп, ограниченной экспоненты. В этой заметке мы исследуем свойства циклически регулярных  полугрупп и связанных с ними многообразий.

\section{Определения и обозначения}

Пусть $X$ - некоторое непустое множетво.  Его элементы будем называть \emph{буквами}, а само множество - \emph{алфавитом}. \emph{Словом} над $X$ назовем любую конечную последовательность букв. Через $S_{X}$ обозначим полугруппу всех слов над $X$ с операцией приписывания.

Слово $u$ называется \emph{повторным}, если каждая буква в нем встречается более одного раза; слово $u$ называется \emph{сингулярным}, если каждая буква встречается в нем только один раз.

Подслово $u^{'}$ слова $u$ называется \emph{циклическим} (или \emph{циклом}), если $\left|{u^{'}}\right|\ge 2$ и $l_{1}\left({u^{'}}\right)\equiv r_{1}\left({u^{'}}\right).$   Слово $u$ называется \emph{покрытым циклами}, если каждая его буква входит в некоторое циклическое подслово слова $u$ . Ясно, что каждое повторное слово покрыто циклами и произведение слов $u\equiv u_{1}u_{2}\cdots u_{m},$ покрытых циклами, само покрыто циклами. Напомним, что через $l_{k}\left({w}\right)$ ($r_{k}\left({w}\right)$) обозначается $k$-я слева (справа) буква слова $w,$ а через $\left|{w}\right|$ обозначается количество букв в слове $w.$

Буква $t$ слова $u$ называется \emph{блокирующей}, если это слово имеет вид $u\equiv u_{1}tu_{2},$ причем $t|u_{1},t|u_{2},u_{1}|u_{2}.$ Запись $v|w$ означает, что слова $v$ и $w$ не имеют общих букв. Знак $\equiv$ используется для обозначения совпадающих слов, а знак  $=$ - для обозначения равных элементов в полугруппах. Отметим, что слово не содержит блокирующих букв тогда и только тогда, когда оно покрывается циклами.

На последовательности букв любого слова $u\equiv a_{i_{1}}a_{i_{2}}\cdots a_{i_{n}}$ можно определить отношение $\sim$ такое, что: $a_{i_{k}}\sim a_{i_{l}},$ если буквы $a_{i_{k}},a_{i_{l}}$ входят в некоторый подцикл  слова $u$ или $k=l.$ Транзитивное замыкание $\sim^{t}$ является эквивалентностью, которая разбивает слово $u$ на подслова $u\equiv u_{1}u_{2}\cdots u_{m},$ состоящие из эквивалентных букв. Ясно, что $u_{i}|u_{j}$ для любых $i\ne j$ .  Такое представление слова будем называть \emph{общим каноническим разложением.}

Количество компонент $m$ в каноническом  разложении слова $u$ будем называть \emph{общей циклической характеристикой} слова и обозначать $m_{c}\left({u}\right).$ Если $m_{c}\left({u}\right)=1$ и $\left|{u}\right|>1,$  то слово $u$ будем называть \emph{регулярным.} Легко видеть, что слово $u$ регулярно тогда и только тогда, когда его длина больше 1 и каждое его двухбуквенное подслово содержится в некотором подцикле.

Ясно, что слово покрывается  циклами тогда и только тогда, когда $\left|{u_{i}}\right|\ge 2$ для всех $i.$ Также отметим, что если $\left|{u_{i}}\right|>1$ для некоторого $i,$ то $u_{i}$ - регулярное слово.

Если какая нибудь буква $x$ встречается в слове $u$ несколько раз, то все ее вхождения попадают в одну компоненту $u_{i}$ ее каконического разложения. Номер $i$ этой компоненты будем называть циклическим номером буквы $x$ в слове $u$ и обозначать $N_{u}\left({x}\right)$ .

На множестве всех букв слова $u$ определяются два бинарных отношения $e_{u}$ и $\le_{u}$ определённые следующим образом:

\begin{align}\nonumber
  xe_{u}y \text{ тогда и только тогда, когда }N_{u}\left({x}\right)=N_{u}\left({y}\right)\\ \nonumber
  x\le _{u}y \text{ тогда и только тогда, когда } N_{u}\left({x}\right)\le N_{u}\left({y}\right).
\end{align}

Два слова  $u,v$ называются \emph{гомогенными}, если они состоят из одних и тех же букв. Два слова $u\equiv u_{1}u_{2}\cdots u_{m_{c}\left({u}\right)}, v\equiv v_{1}v_{2}\cdots v_{m_{c }\left({v}\right)}$ называются \emph{подобными}, если выполнены следующие три условия:

\begin{itemize}
  \item циклические характеристики слов $u,v$ равны, то есть $m_{c}\left({u}\right)=m_{c}\left({v}\right);$
  \item соответствующие компоненты $u_{i}$ и $v_{i}$ являются гомогенными;
  \item если  $\left|{u_{i}}\right|=1,$ то и $\left|{v_{i}}\right|=1$ (и наоборот).
\end{itemize}
Нетрудно видеть, что, если слова $u,v$ подобны, то $e_{u}=e_{v}$ и $\le_{u}=\le_{v}.$

Мы также будем пользоваться следующей четырехэлементной полугруппой, заданной определяющими соотношениями следующим образом:
$$A_{0}=\langle a,b \mid a^{2}=a, b^{2}=b, ba=0\rangle.$$
Если $S$  - некоторая полугруппа, то через $rS$ обозначается множество регулярных элементов $S.$

Полугруппа $S$ называется \emph{регулярно замкнутой,} если произведение регулярных элементов является регулярным элементом, другими словами, если $rS$ - подполугруппа  $S.$ Стоит отметить, что не каждая полугруппа является регулярно замкнутой (например, $A_{0}$). Но тем не менее $A_{0}$ - циклически регулярная полугруппа (но не регулярная).

Однако, известно, если $S$ - ортодоксальная (то есть произведение любых идемпотентов  $S$ являтся идемпотентом), то $S$ регулярно замкнута. Действительно, пусть $a=axa, b=byb$ для некоторых элементов $a,b,x,y$ из $S.$ Тогда $xa, by$ - идемпотенты. Поэтому $xa\cdot by$ - тоже идемпотент. Получим  $a\cdot b=axa\cdot byb=a\left({xa\cdot by}\right)b=a\left({xa\cdot by}\right)\left({xa\cdot by}\right)b=$ $=\left({axa}\right)\left({byxa}\right)\left({byb}\right),$ откуда $ab=abyxab$

Если $S$ -  некоторая полугруппа, то через $I\left({a}\right)$ обозначается главный идеал, порожденный элементом $a,$ то есть $I\left({a}\right)=S^{1}aS^{1}.$

\section{Вспомогательные утверждения}

\begin{lemma}\label{l1}
  Если в полугруппе $S$ для некоторого натурального $n$ выполняются тождества
  \begin{align}\label{eq1}
    x^{2}=x^{n+2}\\ 
    \label{eq2}xyx=(xy)^{n+1}x
  \end{align}
  то $S$ - циклически регулярная полугруппа.
\end{lemma}
\begin{proof}
  Вытекает из определений.
\end{proof}

\begin{lemma}\label{l2}
  Пусть $S$ - циклически регулярная полугруппа, $a,e,f,w$  - элементы $S.$ Если $e,f$  - идемпотенты и принадлежат главному идеалу $I\left({a}\right),$ а также имеет место равенство $w=ew=wf,$ то $w$ - регулярный элемент.
\end{lemma}
\begin{proof}
  По условию идемпотенты $e,f$ принадлежат $I\left({a}\right).$ Поэтому $e=e_{1}ae_{2},$ $f=f_{1}af_{2}.$ Без ограничения общности можно считать $e_{1,2},f_{1,2}$ непустыми символами (в противном случае равенства $e=e_{1}ae_{2},$ $f=f_{1}af_{2}$ можно слева и справа домножить соответственно на $e$ и $f$). Получим, $w=ewf=e_{1}ae_{2}\cdot w\cdot f_{1}af_{2}=e_{1}\left({ae_{2}w\cdot f_{1}a}\right)z\left({ae_{2}w\cdot f_{1}a}\right)f_{2}=ew\left({f_{1}azae_{2}}\right)wf,$ откуда $w=w\left({f_{1}azae_{2}}\right)w.$
\end{proof}

\begin{proposition}\label{p1}
  Главные иделы в циклически регулярной полугруппе регулярно замкнуты.
\end{proposition}
\begin{proof}
  Пусть $S$  - циклически регулярная полугруппа, $a$  - элемент $S.$ Пусть $b,c$  - регулярные элементы $I\left({a}\right).$ Докажем, что $bc$ - регулярный элемент в полугруппе $I\left({a}\right).$ По условию $b=bxb, c=cyc$ для некоторых $x,y$ из $I\left({a}\right).$ Пусть $e=bx, f=yc$ . Ясно, что $e,f$  - идемпотенты и принадлежат $I\left({a}\right).$ С другой стороны,  $e\cdot bc=bxbc=bc,$  $bc\cdot f=bcyc=bc$. По лемме \ref{l2} $w=bc$ - регулярный элемент $S.$ Это означает, что $bc=bczbc$ для некоторого $z$ из  $S.$ Тогда $bc=bc\cdot z^{\#}bc$ , где $z^{\#}=zbcz$ принадлежит $I\left({a}\right)$ Это означает, что $bc$ - регулярный элемент в $I\left({a}\right).$
\end{proof}

\begin{proposition}\label{p2}
  Пусть $S$ - произвольная циклически регулярная полугруппа. Если $w$ - произвольное регулярное слово, то $w$ - регулярный элемент $S$ (имеется ввиду значение $w$ при любой подстановке значений букв из $S$).
\end{proposition}
\begin{proof}
  В дальнейшем слово будем называть \emph{хорошим}, если его значение при любой подстановке значений букв из $S$ является регулярным элементом. Требуется доказать, что любое регулярное слово $w$ является хорошим. Так как $w$ - регулярное слово, то $\left|{w}\right|>1$ по определению. Пусть $x=l_{1}\left({w}\right),$  $y=l_{2}\left({w}\right).$ Так как $w$  -  регулярное слово, то первая и вторая буквы этого слова должны войти в некоторый цикл $w^{'}\le w,$ тогда $u\equiv w'w'',$ где  $w''$ - некоторое слово или пустой символ. По условию,  $w^{'}$ - хорошее слово.

  Пусть слово  $w$ представимо в виде  $w\equiv u\cdot v$ , где $u$ - хорошее слово, а $v$ - слово или пустой символ. Причем,  слово $u$ имеет наибольшую возможную длину. Если $w\equiv u,$ то $w$  - регулярный элемент $S,$ то есть $w$  - хорошее слово.

  Пусть $w\not\equiv u.$ Пусть $y$ - последняя буква в слове $u;$  $z$  - первая буква в слове $v.$ Если слова $u,v$ не имеют общих букв, то $y$ и $z$ не связаны отношением $\sim^{t}$ в слове  $w$ , то есть $w$  - не регулярное слово (так как  нет ни одного цикла с началом в $u$ и концом в $v$). Значит, слова $u$ и $v$ имеют общую букву. Пусть $z$  - общая буква слов $u$  и $v.$ Тогда $u\equiv u_{1}zu_{2}$ ; $v\equiv v_{1}zv_{2}$ , для некоторых слов $u_{1,2}$ и $v_{1,2}.$ Рассмотрим слово $u^{\#}\equiv u_{1}zu_{2}v_{1}z$ . Покажем, что $u^{\#}$ - хорошее слово. Зафиксируем какую-нибудь подстановку букв слова $w$ в полугруппу $S.$ По предположению,  значение слова $u$  - регулярный элемент в полугруппе $S$, то есть в этой полугруппе $u=uau$ (для некоторого $a$ , принадлежащего $S$). Здесь и далее значение любого слова (подслова слова $w$ ) при фиксированной выше подстановке букв и само слово будут обозначаться одинаково, а знак $=$ будет использоваться для обозначения равных элементов в полугруппе $S.$

  Обозначим в полугруппе $S$ элемент $ua$ через $e$ , то есть  $e=ua=u_{1}zu_{2}a.$ Тогда $eu^{\#}=eu_{1}zu_{2}v_{1}z=euv_{1}z=uauv_{1}z=uv_{1}z=u_{1}zu_{2}v_{1}z=u^{\#}.$ По условию, $zu_{2}v_{1}z$ - регулярный элемент в $S,$ то есть  $zu_{2}v_{1}z=\left({zu_{2}v_{1}z}\right)b\left({zu_{2}v_{1}z}\right),$ для некоторого элемента $b$ из $S.$ Обозначим  через $f=bzu_{2}vz_{1}.$ Тогда:
  $$u^{\#}f=\left({u_{1}zu_{2}v_{1}z}\right)bzu_{2}vz_{1}=u_{1}\left({zu_{2}v_{1}z}\right)b\left({zu_{2}v_{1}z}\right)=$$
  $$=u_{1}\left({zu_{2}v_{1}z}\right)=u^{\#}.$$
  Так как элементы $e$ и $f$ имеют в полугруппе $S$ общий делитель $z$ , то они принадлежат  идеалу $I\left({z}\right).$ По Лемме \ref{l2},  $u^{\#}$  - регулярный элемент $S.$ Следовательно, слово $u^{\#}$  - хорошее. Но, так как $w=uv=u_{1}zu_{2}v_{1}zv_{2}=u^{\#}v_{2}$ и $\left|{u^{\#}}\right|>\left|{u}\right|,$ то получаем противоречие с предположением, что длина слова $u$ - наибольшая из возможных. Следовательно, случай  $w\not\equiv u$ невозможен.
\end{proof}

\begin{corollary}\label{c1}
  Значения регулярных слов в полугруппе являются регулярными элементами тогда и только тогда, когда эта полугруппа циклически регулярна.
\end{corollary}

\begin{lemma}\label{l3}
  Если $u$ - регулярное слово, то при любом эндоморфизме $f:S_{X}\rightarrow S_{X}$ образ $u^{f}$ - регулярное слово (здесь и далее образ  $f\left({u}\right)$ слова $u$ обозначается через  $u^{f}$ ).
\end{lemma}
\begin{proof}
  Достаточно доказать, что любые две смежные буквы в слове $u^{f}$ принадлежат некоторому его подциклу. Итак, пусть $xy$  - двухбуквенное подслово $u^{f}.$ Возможны два случая:
  \begin{enumerate}
    \item $xy$  - подслово слова $c^{f}$ для некоторой буквы $c$ слова $u$. Поскольку регулярное слово покрывается циклами, то буква $c$ входит в некоторый подцикл $u^{'}$ слова $u$ .
    \item $x=r_{1}\left({a^{f}}\right), y=l_{1}\left({b^{f}}\right)$ для некоторого $ab$  - двухбуквенного подслова $u.$ Тогда слово $ab$ входит в некоторый подцикл $u^{'}$ слова $u.$
  \end{enumerate}
  Итак, в любом из возможных случаев мы приходим к некоторому подслову $u^{'}$ слова $u$ , являющемуся циклом, причем двухбуквенное слово $xy$ является подсловом слова $\left({u'}\right)^{f}.$ Пусть $u'=xu^{\#}x,$ $x^{f}=a_{i_{1}}a_{i_{2}}\cdots a_{i_{n}},$ тогда $\left({u^{'}}\right)^{f}=a_{i_{1}}a_{i_{2}}\cdots a_{i_{n}}\left({u^{\#}}\right)^{f}a_{i_{1}}a_{i_{2}}\cdots a_{i_{n}}.$ Легко видеть, что любые две смежные буквы слова $\left({u^{'}}\right)^{f}$ входят в некоторый подцикл. В частности, это верно и для букв $x,y.$
  Итак, любые две смежные буквы слова $u^{f}$ принадлежат некоторому подциклу. Следовательно, $u^{f}$ - регулярное слово.
\end{proof}

\begin{proposition}\label{p3}
  Если слова $u$ и $v$  подобны, то тождество $u=v$ выполняеся в полугруппе $A_{0}.$
\end{proposition}
\begin{proof}
  Предположим, что слова $u$ и $v$ подобны. Зафиксируем какую-нибудь подстановку  $f$ из алфавита $X$ в полугруппу $A_{0}.$ По определению слова $u,v$ должны иметь вид: $u=u_{1}u_{2}\cdots u_{m},$ $v=v_{1}v_{2}\cdots v_{m},$ причем,  соответствующие компоненты $u_{i}$ и $v_{i}$ являются гомогенными. Также, на основании определения, если $\left|{u_{i}}\right|=1,$ то и $\left|{v_{i}}\right|=1$ (и наоборот) для всех $i=1,2,\cdots,m.$ Рассмотрим возможные случаи.
  \begin{enumerate}
    \item $\left|{u_{i}}\right|>1.$ В таком случае $\left|{v_{i}}\right|>1.$ Слова  $u_{i},v_{i}$ являются регулярными (по определению канонического разложения слова). По предложению \ref{l2}, в любой циклически регулярной полугруппе значения регулярных слов являются регулярными элементами. Легко заметить, что полугруппа $A_{0}$  - циклически регулярная. Действительно, эта полугруппа состоит из четырех элементов $a,b,z=ab$ и $0$ и может быть задана таблицей Кэли: $A_{0}=\{a,b,z,0\mid a^{2}=a, b^{2}=b, ab=az=zb=z\}.$ Неуказанные произведения считаются равными нулю. Нетрудно убедиться, что в этой полугруппе выполняются тождества \eqref{eq1} и \eqref{eq2} для $n=1.$ По лемме \ref{l1} заключаем, что $A_{0}$  - циклически регулярная полугруппа. Поэтому значения слов  $u_{i}$ и $v_{i}$ должны быть регулярными элементами в полугруппе  $A_{0}.$ Следует отметить, что, если при рассматриваемой фиксированной подстановке значение хотя бы одной буквы $d$  из $u_{i}$ равно $z$ , то  значение всего слова $u_{i}$ обязано равняться $0$ .  Действительно, поскольку $u_{i}$  - регулярное слово, то оно должно иметь вид  $u\equiv w_{1}cw_{2}dw_{3}cw_{4}$ для некоторой буквы $c$ и некоторых слов $w_{1,2,3.4}$ (возможно пустых). Тогда $f\left({u}\right)=f\left({w_{1}}\right)f\left({c}\right)f\left({w_{2}}\right)zf\left({w_{3}}\right)f\left({c}\right)f\left({w_{4}}\right)=0.$ Это же утверждение имеет место и для слова $w_{i}.$ Регулярными элементами в полугруппе $A_{0}$ являются только $a, b, 0.$ Если значение слова $u_{i}$ равно $a,$ то значение всех букв слова $u_{i}$ равно $a.$ Тогда значение всех букв слова $v_{i}$ тоже равно $a,$ поскольку слова $u_{i}, v_{i}$ состоят из одинаковых букв. Если значение слова $u_{i}$ равно $b,$ то значение всех букв слова $u_{i}$ равно $b.$ Тогда значение всех букв слова $v_{i}$ тоже равно $b,$ поскольку слова $u_{i},v_{i}$ состоят из одинаковых букв. Если значение слова $u_{i}$ равно $b,$ то значение всех букв слова $u_{i}$ равно $b.$ Тогда значение всех букв слова $v_{i}$ тоже равно $b,$ поскольку слова $u_{i},v_{i}$ состоят из одинаковых букв. Если значение слова $u_{i}$ равно $0$ , то возможно два случая:
        \begin{enumerate}
                   \item\label{sl1a} значение одной из букв слова $u_{i}$ равно $0;$
                   \item\label{sl1b} значение ни одной из букв слова $u_{i}$ не равно $0.$
        \end{enumerate}
        Если имеет место случай \ref{sl1a}, то значение слова $v_{i}$ равно $0,$ поскольку слова $u_{i}, v_{i}$ состоят из одинаковых букв.

        Пусть имеет место случай \ref{sl1b}. Если значение хотя бы из одной букв слова $u_{i}$ равно $z,$  то, как отмечалось выше, значения слов $u_{i}$ и $v_{i}$ равны $0,$ поскольку эти слова состоят из одинаковых букв. Осталось рассмотреть ситуацию, при которой значения букв слов $u_{i}$  равны либо $a,$  либо $b,$ причем в состав слова $f\left({u_{i}}\right)=u_{i}\left({a,b}\right)$ обязательно входят буквы $a $ и $b.$ Тогда  значение слова $v_{i}$ равно $f\left({v_{i}}\right)=v_{i}\left({a,b}\right),$ причем в состав слова $v_{i}\left({a,b}\right)$ входят буквы $a$ и  $b.$ Поскольку $f\left({v_{i}}\right)=v_{i}\left({a,b}\right)$  - регулярный элемент, то $v_{i}\left({a,b}\right)=v_{i}\left({a,b}\right)v^{\#}v_{i}\left({a,b}\right)$ для некоторого элемента $v^{\#}$ из $A_{0}.$ Если $v^{\#}=0,$ то $v_{i}\left({a,b}\right)=0=u.$ Если $v^{\#}\ne 0,$ то $v^{\#}=v^{\#}\left({a,b}\right)$ - слово от элементов $a,b.$ Заключаем, $v_{i}\left({a,b}\right)=v_{i}\left({a,b}\right)v^{\#}v_{i}\left({a,b}\right)=v_{i}\left({a,b}\right)v^{\#}\left({a,b}\right)v_{i}\left({a,b}\right)$ содержит подслово $ba.$ Это означает, что $v_{i}\left({a,b}\right)=0=u_{i}\left({a,b}\right).$ Итак получаем, что значения слов $u_{i}$ и $v_{i}$ совпадают (при любой подстановке элементов  $A_{0}$ вместо букв). Это означает, что в полугруппе $A_{0}$ выполняется тождество  $u_{i}=v_{i}.$
    \item $\left|{u_{i}}\right|=1.$ Тогда $\left|{v_{i}}\right|=1$ (по определению подобных слов). Поскольку слова $u_{i},v_{i}$ гомогенные,  заключаем,  $u_{i}=v_{i}$ - одинаковые буквы. Итак, во всех возможных случаях в полугруппе $A_{0}$ выполняется тождество $u_{i}=v_{i}.$
  \end{enumerate}
  Откуда следует, что в полугруппе $A_{0}$ выполняется тождество $u=v.$
\end{proof}

\section{Основные результаты}

В следущей теореме используются семь конечных полугрупп из работы автора \cite{1}, заданных определяющими соотношениями:

$A=\langle x,y \mid x=x^{2}; y^{2}=0; xy=yx\rangle,$

$B=\langle x,y \mid x^{2}=0; y^{2}=0; xyx=yxy\rangle,$

$C_{\lambda}=\langle x,y \mid x^{2}=x^{3}; xy=y; yx^{2}=0; y^{2}=0\rangle,$

$C_{\rho}=\langle x,y \mid x^{2}=x^{3}; yx=y; x^{2}y=0; y^{2}=0\rangle,$

$N_{3}=\langle x \mid x^{3}=0\rangle,$

$D=\langle x,y \mid x^{2}=0; y=y^2; yxy=0\rangle,$

$K_{n}=\langle x,y \mid x^{2}=0; y^{2}=y^{n+2}; yxy=0; xy^{q}x=0, (q=2,\cdots,n); xyx=xy^{n+1}x\rangle.$

\begin{theorem}\label{t1}
  Для многообразия полугрупп $V$ следующие войства эквивалентны.
  \begin{enumerate}
    \item Все полугруппы многообразия $V$ цикличски регулярны.
    \item В многообразии $V$ выполняются тождества \eqref{eq1} и \eqref{eq2} для некоторого натурального числа $n.$
    \item $V$ не содержит ни одной из конечных полугрупп $A,$ $B,$ $C_{\lambda},$ $C_{\rho},$ $N_{3},$ $D$ и $K_{n}.$
  \end{enumerate}
\end{theorem}
\begin{proof}
  Следует из леммы \ref{l1} и леммы 10 из работы \cite{1}. Достаточно проверить, что перечисленные выше конечные полугруппы не являются циклически регулярными.
\end{proof}

\begin{corollary}\label{c2}
  Если конечные полугруппы многообразия являются циклически регулярными, то и все полугруппы этого многообразия являются циклически регулярными.
\end{corollary}

\begin{corollary}\label{c3}
  Существует алгорит полиномиальной сложности, определяющий по конечному базису тождеств многообразия, все ли его полугруппы являются циклически регулярными.
\end{corollary}
\begin{proof}
  Следует из теремы \ref{t1} и того факта, что проверка выполнимости тождеств на полугруппах  $\left({1}\right)-\left({7}\right)$ имеет полиномиальную сложность (см. \cite{2}).
\end{proof}

\begin{theorem}\label{t2}
  Для многообразия полугрупп $V$ следующие свойства эквивалентны.
  \begin{enumerate}
    \item\label{t2_1} Все полугруппы многообразия $V$  регулярно замкнуты.
    \item\label{t2_2} Среди тождеств базиса многообразия $V$ есть тождество $u=v,$ у которого слова $u,v$ не подобны.
    \item\label{t2_3} В многообразии $V$ выполняется тождество вида $x^{n}y^{n}=uyxv,$ где $x,y$ -буквы, а $u,v$ - слова (возможно пустые).
    \item\label{t2_4} Многообразие $V$ не содержит полугруппы $A_{0}.$
  \end{enumerate}
\end{theorem}
\begin{proof}
  \ref{t2_1}$\Rightarrow$\ref{t2_2}. Предположим, все полугруппы многообразия $V$ регулярно замкнуты, но все тождества $u=v$ базиза, которым задается $V,$ состоят из подобных слов. В этом случае $A_{0}$ принадлежит многообразию $V$ по предложению \ref{p3}. Но полугруппа $A_{0}$ не является регулярно замкнутой, т.к. произведение регулярных элементов $a\cdot b$ не является регулярным элементом. Противоречит предположению.

  \ref{t2_2}$\Rightarrow$\ref{t2_3}. Пусть среди тождеств базиса многообразия $V$ выполняется тождество $u=v$ , у которого слова $u,v$ не подобны. Это возможно в одном из четырех случаев:
  \begin{enumerate}
    \item\label{sl2.1} $u,v$ - не гомогенные слова;
    \item\label{sl2.2} $u,v$ - гомогенные слова, но $e_{u}\ne e_{v};$
    \item\label{sl2.3} $u,v$ - гомогенные слова, $e_{u}=e_{v},$ но $\le_{u}\ne\le _{v};$
    \item\label{sl2.4} $u,v$ - гомогенные слова, $e_{u}=e_{v},$ $\le_{u}=\le_{v},$ причем, для некоторого $i$ слова  $u_{i}$ и $v_{i}$ имеют вид $u_{i}=z,v_{i}=z^{k}$ для некоторой буквы $z$ и натурального числа  $k>1.$
  \end{enumerate}
  В первом случае \ref{sl2.1}, получаем, что в качестве следствия тождества $u=v$ можно получить тождество вида $u^{\#}yv^{\#}=x^{n}.$ Подставляя вместо $y\mapsto xy$ и умножая справа полученное тождество на $y^{n},$ получим в качестве следствия тождество вида $x^{n}y^{n}=uyxv.$

  Во втором случае \ref{sl2.2} получаем, что существует пара различных букв $x,y$ , для которой $N_{u}\left({x}\right)\ne N_{u}\left({y}\right),$ но $N_{v}\left({x}\right)=N_{v}\left({y}\right)$ (или наоборот). Пусть для определенности $N_{u}\left({x}\right)<N_{u}\left({y}\right).$ Рассмотрим отображение $F:X_{u}\rightarrow X_{u}$ (где $X_{u}$ - множество букв в слове $u$), определенное правилом
  $$f(z)=\begin{cases}
    x &\text{ если } N_{u}\left({z}\right)\le N_{u}\left({x}\right);\\
    y &\text{ если } N_{u}\left({z}\right)>N_{u}\left({x}\right).
  \end{cases}$$
  Заметим, что $f\left({x}\right)=x, f\left({y}\right)=y.$ Тогда $u^{f}=u^{f}_{1}u^{f}_{2}..u^{f}_{m}=x^{k}y^{l}.$ Поскольку $N_{v}\left({x}\right)=N_{v}\left({y}\right),$ то буквы  $x,y$ лежат в одной компоненте $v_{i}$ канонического разложения слова $v=v_{1}..v_{i}..v_{m}.$ По определению $v_{i}$  - регулярное слово. Получим, $v^{f}=\left({v_{1}..}\right)^{f}.v^{f}_{i}\left({..v_{m}}\right)^{f}=v^{'}.v^{f}_{i}\cdot v^{''}.$ Отметим, что $v^{'},v^{f}_{i},v^{''}$  - слова от букв $x,y,$ причем, буквы $x,y$ присутствуют в слове $v^{f}_{i}.$ По лемме \ref{l3}, $v^{f}_{i}$  - регулярное слово (как образ регулярного) от переменых $x,y.$ Заключаем, что слово $yx$  - подслово $v^{f}_{i}.$ В противном случае слово $v^{f}_{i}=x^{p}y^{q}$ не является регулярным. Вывод: следствием тождества $u=v$ является тождество вида $x^{k}y^{l}=u^{\#}yxv^{\#}.$

  В третьем случае \ref{sl2.3}, получаем, что существует пара различных букв $x,y$ , для которой $N_{u}\left({x}\right)<N_{u}\left({y}\right),$ но $N_{v}\left({x}\right)>N_{v}\left({y}\right)$ (или наоборот). Рассмотрим отображение $F:X_{u}\rightarrow X_{u}$ (где $X_{u}$  - множество букв в слове $u$ ), определенное правилом
  $$f(z)=\begin{cases}
    x &\text{ если } N_{u}\left({z}\right)\le N_{u}\left({x}\right);\\
    y &\text{ если } N_{u}\left({z}\right)>N_{u}\left({x}\right).
  \end{cases}$$
  Заметим, что $f\left({x}\right)=x , f\left({y}\right)=y.$ Тогда $u^{f}=u^{f}_{1}u^{f}_{2}..u^{f}_{m}=x^{k}y^{l},$ $v^{f}=v^{'}yv^{''}xv^{'''}.$ Вывод: следствием тождества $u=v$ является тождество вида $x^{k}y^{l}=u^{\#}yxv^{\#}.$

  В четвёртом случае \ref{sl2.4}, рассмотрим отображение $F::X_{u}\rightarrow X_{u}$ (где $X_{u}$ - множество букв в слове $u$ ), определенное правилом
  $$f(t)=\begin{cases}
    x &\text{ если } N_{u}\left({t}\right)\le i\\
    y &\text{ если } N_{u}\left({t}\right)>i\\
    xy &\text{ если } t=z.
  \end{cases}$$
  Тогда $u^{f}=u^{f}_{1}u^{f}_{2}..u^{f}_{m}=x^{k}y^{l},$ $v^{f}=v^{'}yv^{''}xv^{'''}.$ Вывод: и в этом случае следствием тождества $u=v$ является тождество вида  $x^{k}y^{l}=u^{\#}yxv^{\#}.$

  \ref{t2_3}$\Rightarrow$\ref{t2_4}. Достаточно показать, что тождество вида $x^{n}y^{n}=uyxv$ не выполняется в полугруппе  $A_{0}.$  Действительно, при подстановке $x\rightarrow a, y\rightarrow b$ значение левой части равно $ab,$ а значение правой части равно $0.$

  \ref{t2_3}$\Rightarrow$\ref{t2_1}. Пусть в многообразии $V$ выполняется тождество вида  $x^{n}y^{n}=uyxv$ . Пусть $a,b$  - регулярные элементы некоторой полугруппы $S$ из $V.$ Имеем $a=aa^{'}a, b=bb^{'}b$ для некоторых $a^{'},b^{'}$ из $S$. Пусть $e=a^{'}a, f=bb^{'}$ . Ясно, что $e,f$  - идемпотенты полугруппы $S.$ При подстановке $x\rightarrow e ; y\rightarrow f$ в тождесво $x^{n}y^{n}=uyxv$ получим равенство $ef=cfed$ в полугруппе $S$ для некотрых элементов $c,d$ из  $S,$ причем, $c=c\left({e,f}\right),d=d\left({e,f}\right)$ - это выражения от элементов $e,f$ или пустые символы. Откуда $ef=ec\left({e,f}\right)fed\left({e,f}\right)f=efwef$ для некоторого элемента $w$ из $S.$ Далее, $ab=aa^{'}a\cdot bb^{'}b=a\left({a^{'}a\cdot bb^{'}}\right)w\left({a^{'}a\cdot bb^{'}}\right)b=ab\left({b^{'}wa^{'}}\right)ab.$ Следовательно, $ab$  - регулярный элемент.

  \ref{t2_4}$\Rightarrow$\ref{t2_1}. Пусть многообразие $V$ не содежит полугруппы $A_{0}.$ Предположим, что у всех тождеств $u=v$ из базиса, слова $u,v$ - подобные. По предложению \ref{p3}, $A_{0}$ принадлежит многообразию $V.$ Противоречение.
\end{proof}

\begin{corollary}\label{c4}
  Существует алгоритм линейной сложности (от входных данных), определяющий регулярную замкнутость всех полугрупп многообразия, заданного конечным набором тождеств.
\end{corollary}

\begin{corollary}\label{c5}
  Тождество $u=v$ выполняется на полугруппе $A_{0}$ тогда и только тогда ,когда слова $u,v$ подобны.
\end{corollary}

\bibliographystyle{amsplain}

\end{document}